\documentclass[preprint,11pt]{elsarticle}

\newtheorem{theorem}{Theorem}[section]
\newtheorem{corollary}{Corollary}[section]

\usepackage{amssymb}

\begin{document}

\begin{frontmatter}

\title{Self-normalized deviation inequalities   with application to $t$-statistics}
\author{Xiequan Fan$^*$}
 \cortext[cor1]{\noindent Corresponding author. \\
  \mbox{ \ \ \ } \textit{E-mail}: fanxiequan@hotmail.com (X. Fan). }
\address{Center for Applied Mathematics, Tianjin University, 300072 Tianjin,  China}

\begin{abstract}
Let $(\xi _i)_{i=1,...,n}$ be a sequence of independent and symmetric random variables. We consider the upper bounds on tail probabilities
of self-normalized deviations $$ \mathbf{P} \Big( \max_{1\leq k \leq n} \sum_{i=1}^{k} |\xi_i|\big/ \big(\sum_{i=1}^{n} |\xi_i|^\beta \big)^{1/\beta}   \geq x \Big) $$
for $x>0$ and  $\beta >1.$ Our bound is the  best that can be obtained from the  Bernstein inequality under the present assumption. An application to
Student's $t$-statistics is also given.
\end{abstract}

\begin{keyword} Self-normalized deviations; Student's $t$-statistic; exponential inequalities
\vspace{0.3cm}
\MSC Primary 60E15; 60F10
\end{keyword}

\end{frontmatter}


\section{Introduction}
Let $(\xi _i)_{i\geq1}$ be a sequence of independent, centered and nondegenerate real-valued random variables (r.v.s). Denote by $$ S_n=\sum_{i=1}^n \xi_i  \ \ \ \ \textrm{and}\ \ \ \ V_n(\beta) =\Big(\sum_{i=1}^{n} |\xi_i|^\beta \Big)^{1/\beta}, \ \beta>1.$$
 The study of the tail probabilities $\mathbf{P}( S_n/V_n(\beta) \geq x ) $ certainly has attracted some particular attentions.
 In the case where r.v.s\ $(\xi _i)_{i\geq1}$ are identically distributed and $\mathbf{E}|\xi_1|^\beta=  \infty, \beta>1$, Shao \cite{S97} proved
 the following deep large deviation principle (LDP) result: for any $x> 0,$
 $$ \lim_{n\rightarrow \infty}\mathbf{P}\Big( \frac{ S_n}{V_n(\beta)\, n^{1-1/\beta}}  \geq x\Big)^{1/n}
 =\sup_{c\geq0} \inf_{t\geq 0} \mathbf{E}\Big[ \exp\Big\{t\Big(cX-x\Big(\frac{1}{\beta }|X|^\beta + \frac{\beta-1}{\beta}c^{\beta/(\beta-1)} \Big ) \Big) \Big\} \Big].$$
 The related moderate deviation principles (MDP) are also given by Shao \cite{S97} and Jing,   Liang and Zhou \cite{JLZ12}.
However,  the LDP and MDP results  do not diminish the need for tail probability inequalities valid for given $n.$ Such inequality have been obtained in particular
by  Wang and Jing \cite{WJ99}.   They proved that if the r.v.s\ $(\xi _i)_{i\geq1}$ are symmetric (around $0$), then for all $x >0,$
\begin{eqnarray} \label{f12}
 \mathbf{P}\bigg(   \frac{S_n}{V_n(2)}  \geq x \bigg)
  \leq    \exp\left\{- \frac{x^2}{ 2  }   \right\}.
\end{eqnarray}
This bound is rather tight for moderate $x$'s. Indeed, as showed by the MDP  result of  Shao \cite{S97} (cf.\ Theorem 3.1), for certain class of r.v.s  it holds that
\begin{eqnarray}
 x_n^{- 2   }\ln  \mathbf{P}\bigg(   \frac{S_n}{V_n(2)}  \geq x_n \bigg)
 =   -\frac12,
\end{eqnarray}
with $x_n \rightarrow \infty$ and $x_n =o(\sqrt{n}).$  See also Theorem 2.1 of Jing,   Liang and Zhou \cite{JLZ12} for non identically distributed r.v.s.
In Fan, Grama and Liu   \cite{F15}, inequality (\ref{f12}) has been further extended to the case of partial maximum: for all $x>0,$
\begin{eqnarray}\label{f13}
 \mathbf{P}\bigg( \max_{1\leq k \leq n} \frac{S_k}{ V_n(2)}  \geq x\bigg)
   \leq  \exp \left\{-  \frac{x^2}{2 \, }    \right\}.
\end{eqnarray}
On the other hand, by the Cauchy-Schwarz inequality, it is easy to see that  $S_n^2 \leq n\, (V_n(2))^2.$ Thus for all $x >\sqrt{n},$
\begin{eqnarray}\label{f14}
\mathbf{P}\bigg(   \frac{S_n}{ V_n(2)} \geq x \bigg) \leq  \mathbf{P}\bigg(   \max_{1\leq k \leq n} \frac{S_k}{ V_n(2)} \geq x \bigg)  =  0,
\end{eqnarray}
which cannot be deduced from (\ref{f12}) and (\ref{f13}).
Hence, the inequalities (\ref{f12}) and (\ref{f13}) are not tight enough.

In this paper we  give an improvement on inequality (\ref{f13}). Our inequality  implies (\ref{f14}).  More general, we establish
  the upper bound on tail probabilities $\mathbf{P}\big( \max_{1\leq k \leq n} S_k/ V_n(\beta)   \geq x \big), x>0,$ for symmetric r.v.s\ $(\xi _i)_{i\geq1}.$
In particular, we show that our
inequality is the best   that can be obtained from the classical Bernstein inequality: $\mathbf{P}(X>x)\leq \inf_{\lambda >0}\mathbf{E}[e^{\lambda(X-x)}].$
An application to Student's $t$-statistics is also given.

The paper is organized as follows. Our main result and the applications are stated and discussed in Section \ref{sec2}.
Proofs are deferred to Section \ref{sec4}.

\section{Main results}  \label{sec2}

In the following theorem, we give a self-normalized deviation  inequality   for  independent and symmetric random variables.
\begin{theorem}\label{th1}
Assume  that $(\xi _i)_{i=1,...,n}$ is a sequence of independent, symmetric and nondegenerate random variables. Given a constant $\beta \in (1, \infty),$ denote by $$V_n(\beta) =\Big(\sum_{i=1}^{n} |\xi_i|^\beta \Big)^{1/\beta}.$$
Then for all $ 0 < x  \leq n^{\frac{\beta-1}{\beta}  }$,
\begin{eqnarray}
 \mathbf{P}\left( \max_{1\leq k \leq n} \frac{S_k}{ V_n(\beta)}    \geq x\right)   \leq   B_n(\beta, x):= \frac{1}{2^n} \bigg( \sqrt{t} + \frac{1}{\sqrt{t}}   \bigg)^n t^{- \frac 12 n^{1/\beta} x } , \label{ineqa1}
\end{eqnarray}
where $$t= \frac{n^{\frac{\beta-1}{\beta}  }+ x }{n^{\frac{\beta-1}{\beta}  }- x}$$
 with the convention  $B_n(\beta, n^{\frac{\beta-1}{\beta}  })=2^{-n}$.
 Moreover  $B_n(2, x)$ is increasing in $n$ and for any $x>0,$
$$\lim_{n\rightarrow \infty}  B_n(2, x) =\exp \left\{-  \frac{x^2}{2 \, }   \right\} . $$
\end{theorem}

The range of validity of our inequality (\ref{ineqa1})  is optimal.  Indeed, H\"{o}lder's inequality implies that
$S_n  \leq V_n(\beta)n^{\frac{\beta-1}{\beta}  }.$
Thus when $x> n^{\frac{\beta-1}{\beta}  },$ it holds that $$\mathbf{P}\Big( \max_{1\leq k \leq n} \frac{S_k}{ V_n(\beta)}  \geq x\Big) =0.$$

Notice that bound (\ref{ineqa1}) is the best that can be obtained from the following Bernstein inequality
\begin{eqnarray}
 \mathbf{P}\bigg(  \frac{S_n}{ V_n(\beta)}  \geq x \bigg)
 \leq \,  \inf_{\lambda\geq 0} \mathbf{E} \Big[ e^{\lambda \big(\frac{S_n}{ V_n(\beta)}- x \big) } \Big].
\end{eqnarray}
Indeed, if $\xi _i=\pm a, a>0,$ with probabilities $1/2,$ then it holds  for all $ 0< x < \sqrt{n}$,
 $$  \inf_{\lambda\geq 0} \mathbf{E} \Big[ e^{\lambda \big(\frac{S_n}{ V_n(\beta)}- x \big) } \Big]  =  \inf_{\lambda\geq 0} \mathbf{E} \Big[ e^{\lambda \big(\frac{S_n}{a  n^{1/\beta}}- x \big) } \Big] =\inf_{\lambda\geq 0}e^{ - \lambda x  } \Big( \cosh(\frac{\lambda}{n^{1/\beta}})    \Big)  ^n = B_n(\beta, x).$$
Moreover,  when $x\rightarrow n^{\frac{\beta-1}{\beta}  },$  bound (\ref{ineqa1}) tends to $2^{-n},$ which
is the best possible at $x=n^{\frac{\beta-1}{\beta}  }.$  Indeed, for the $\xi _i$'s mentioned above,    it holds
 $$ \mathbf{P}\bigg( \max_{1\leq k \leq n} \frac{S_k}{ V_n(\beta)}  \geq n^{\frac{\beta-1}{\beta}  }\bigg) = \mathbf{P}\Big(\xi _i=  a \  \textrm{for all}\ i \in [1, n] \Big)= \frac{1}{2^{n}} .$$

Since the random variables $(\xi _i)_{i=1,...,n}$ are  symmetric,  it  is obvious that  for all $ 0 < x  \leq n^{\frac{\beta-1}{\beta}  }$,
\begin{eqnarray*}
 \mathbf{P}\left( \max_{1\leq k \leq n} \frac{S_k}{ V_n(\beta)}   \leq- x\right)   \leq   B_n(\beta, x),
\end{eqnarray*}
where $B_n(\beta, x)$ is defined by (\ref{ineqa1}).

When $\beta \in (1, 2],$ inequality (\ref{ineqa1}) implies the following bound.

\begin{corollary}\label{co1} Assume the condition of Theorem \ref{th1}.
If $\beta \in (1, 2],$ then for all $x>0,$
\begin{eqnarray}
 \mathbf{P}\bigg( \max_{1\leq k \leq n} \frac{S_k}{ V_n(\beta)}  \geq x\bigg)
   \leq  \exp \left\{-  \frac{x^2}{2 \, } n^{\frac{2}{\beta}- 1   }   \right\}. \label{ineqa2}
\end{eqnarray}
In particular, the last inequality implies that  for $\beta \in (1, 2),$ $$\frac{S_n}{ V_n(\beta)}  \rightarrow 0, \ \ \ \ \ n\rightarrow \infty,$$ in probability.
\end{corollary}

For $\beta \in (1, 2],$ inequality (\ref{ineqa2}) implies the following upper bound of LDP:
 \begin{eqnarray}\label{fds01}
 \limsup_{n \rightarrow \infty}\frac1n \ln \mathbf{P}\bigg( \max_{1\leq k \leq n} \frac{S_k}{ V_n(\beta) n^{\frac{\beta-1}{\beta}  }}  \geq x\bigg)
   \leq  -  \frac{x^2}{2 \, }  , \ \ \ \ \ x \in (0, 1].
\end{eqnarray}
It also implies the following upper bound of MDP: for any $\alpha \in ( \frac{\beta-2}{2\beta }  ,  \frac{\beta-1}{ \beta }   ) $,
 \begin{eqnarray}
 \limsup_{n \rightarrow \infty}\frac1{n^{2\alpha+ \frac2\beta -1}} \ln \mathbf{P}\left( \max_{1\leq k \leq n} \frac{S_k}{ V_n(\beta) n^{\alpha }}  \geq x\right)
   \leq  -  \frac{x^2}{2 \, }  , \ \ \ \ \ x \in (0, \infty).
\end{eqnarray}
With certain regularity conditions on tail probabilities of $\xi_i,$  the LDP and MDP results are allowed to be established. We refer to Shao \cite{S97} and Jing,  Liang and Zhou \cite{JLZ12}.

 Wang and Jing \cite{WJ99}  proved that for all $x >0,$
\begin{eqnarray}\label{fdgfgfdgh1}
 \mathbf{P}\bigg(   \frac{S_n}{V_n(2)}  \geq x \bigg)
  \leq    \exp\left\{- \frac{x^2}{ 2  }   \right\}.
\end{eqnarray}
An earlier result similar to (\ref{fdgfgfdgh1}) can be found in \cite{H90}, where Hitczenko  has obtained the same upper bound    on tail probabilities $\mathbf{P} ( S_n   \geq x \, ||\sqrt{[S]_n} ||_{\infty} ).$  When $\beta=2,$ inequality (\ref{ineqa2}) reduces to the following inequality of Fan \emph{et al.}\ \cite{F15}:
 for all $x>0,$
\begin{eqnarray}\label{fdgfgfdgh2}
 \mathbf{P}\bigg( \max_{1\leq k \leq n} \frac{S_k}{ V_n(2)}  \geq x\bigg)
   \leq  \exp \left\{-  \frac{x^2}{2 \, }    \right\}.
\end{eqnarray}
Thus inequality (\ref{ineqa2})  can be regarded as a generalization of  (\ref{fdgfgfdgh1}) and (\ref{fdgfgfdgh2}).
Moreover, since the bound (\ref{ineqa1}) is less  than the bound (\ref{fdgfgfdgh2}), our inequality (\ref{ineqa1})   improves
 on (\ref{fdgfgfdgh2}).

Let $\{Y_i\}_{ i\geq 1}$ be a sequence of independent nondegenerate r.v.s, and $\{d_i\}_{ i\geq 1}$
be a sequence of independent
Rademacher r.v.s, i.e. $\textbf{P}(d_i=\pm 1)=\frac12.$ Let $\xi_i=d_iY_i.$ Assume that $\{Y_i\}_{ i\geq 1}$ and $\{d_i\}_{ i\geq 1}$
are independent. Then we now  have
$$S_n=\sum_{i=1}^nd_iY_i, \ \ \ \ \ \ V_{n}(\beta)=\Big(\sum_{i=1}^n|Y_i|^\beta\Big)^{1/\beta}, \  \ \ \ \textrm{for} \ \beta> 1.$$
The following result easily follows from   Theorem \ref{th1} and Corollary \ref{co1}.
\begin{corollary} Let $\xi_i=d_iY_i$ for $i=1,...,n$.
If $\beta \in (1, 2],$ then for all $ 0 < x  \leq n^{\frac{\beta-1}{\beta}  },$
\begin{eqnarray}
 \mathbf{P}\left( \max_{1\leq k \leq n} \frac{S_k}{ V_n(\beta)}    \geq x\right) \  \leq  \ B_n(\beta, x) \ \leq \ \exp \left\{-  \frac{x^2}{2 \, } n^{\frac{2}{\beta}- 1   }   \right\} .
\end{eqnarray}
In particular, the last inequality implies that  for $\beta \in (1, 2),$ $$\frac{S_n}{ V_n(\beta)}  \rightarrow 0, \ \ \ \ \ n\rightarrow \infty,$$ in probability.
\end{corollary}

  Consider Student's $t$-statistic $T_n$ defined by
\[
T_n = \sqrt{n} \, \overline{\xi}_n / \widehat{\sigma},
\]
where $$\overline{\xi}_n = \frac{S_n}{n}  \ \ \ \textrm{and}\  \ \ \widehat{\sigma}^2 = \sum_{i=1}^n  \frac{(\xi_i - \overline{\xi}_n )^2}{ n-1}  .$$ It is known that for all $x>0,$
\[
\mathbf{P}\Big( T_n  \geq x \Big) = \mathbf{P}\bigg(  \frac{S_n }{\sqrt{[S]_n}}  \geq x \Big(\frac{n}{n+x^2-1} \Big)^{1/2}  \bigg );
\]
 see Efron \cite{E69}. Notice that for all  $x>0,$ it holds
$0< x \Big(\frac{n}{n+x^2-1} \Big)^{1/2} \leq n^{1/2}.$ With the help of (\ref{ineqa1}), we have the following exponential bound for Student's $t$-statistics.
\begin{theorem}  Assume  that $(\xi _i)_{i=1,...,n}$ is a sequence of independent, symmetric and nondegenerate random variables.
Then  for all $  x> 0$,
\begin{eqnarray}
 \mathbf{P}\Big( T_n  \geq x \Big)
\ \leq \   B_n \bigg(2,\  x \Big(\frac{n}{n+x^2-1} \Big)^{1/2}  \bigg)  ,
\end{eqnarray}
where $B_n(2, x)$ is defined by (\ref{ineqa1}).
\end{theorem}

\section{Proofs of Theorems}\label{sec4}
The proof  of Theorem \ref{th1} is based on a  method  called change of probability measure for martingales.
The method is developed by Grama  and Haeusler \cite{GH00}.

\noindent\emph{Proof of Theorem \ref{th1}.}   For any $i=1,...,n$, set
\begin{eqnarray}\label{dun87}
\eta_i= \frac{ \xi_i}{ V_n(\beta)}, \ \ \mathcal{F}_{0} = \sigma \Big( |\xi_j|, 1\leq j\leq n \Big) \ \  \textrm{and}  \ \ \mathcal{F}_{i} = \sigma \Big( \xi_{k},  1\leq k\leq i,\  |\xi_j|, 1\leq j\leq n \Big).
\end{eqnarray}
Since $(\xi _i)_{i=1,...,n}$ are  independent and  symmetric, then  $$ \mathbf{E}[\xi_i> y\, |\, \mathcal{F}_{i-1} ] = \mathbf{E}\Big[\xi_i> y \, \Big|\, |\xi_i|  \Big] =\mathbf{E}\Big[ -\xi_i > y \,\Big|\, |-\xi_i| \Big]=\mathbf{E}[ -\xi_i > y\, |\, \mathcal{F}_{i-1}  ].$$
Thus $\big(\eta_i,\mathcal{F}_i\big)_{i=1,...,n}$ is a sequence of conditionally symmetric martingale differences, i.e.\ $ \mathbf{E}[\eta_i> y|\, \mathcal{F}_{i-1} ] =  \mathbf{E}[ -\eta_i > y|\, \mathcal{F}_{i-1}  ].$  It is easy to see that
\begin{equation}
\frac{S_n}{\textrm{V}_n(\beta)} = \sum_{i=1}^n \eta_i   \nonumber
\end{equation}
is a sum of martingale differences,  and that $\big(\eta_i,\mathcal{F}_i\big)_{i=1,...,n}$ satisfies  $$\sum_{i=1}^n  |\eta_i|^\beta  \  = \ \sum_{i=1}^n\frac{ |\xi_i|^\beta}{ V_n(\beta)^\beta}  \ =\ 1.$$
For   any $x >0$, define the stopping time $T$:
\[
T(x) =\min\Big\{k\in [1, n]: \sum_{i=1}^k\eta_i \geq x\Big\},
\]
with the convention that $\min{\emptyset}=0$. Then it follows that
 \[
\textbf{1}_{\big\{  \max_{1\leq k \leq n} S_k/ V_n(\beta) \, \geq x    \big\}} = \sum_{k=1}^{n}  \textbf{1}_{\{ T(x) =k\}}.
\]
For any
nonnegative number $\lambda$, define the  martingale $M(\lambda )=(M_k(\lambda ),\mathcal{F}_k)_{k=0,...,n},$ where
\[
M_k(\lambda )=\prod_{i=1}^k\frac{\exp\left\{\lambda \eta_i    \right\}}{\mathbf{E}\left[\exp\left\{\lambda \eta_i   \right\} |
\mathcal{F}_{i-1} \right]},\quad \quad  \quad M_0(\lambda )=1.
\]
Since $T$ is a stopping time, then $M_{T\wedge k}(\lambda )$, $\lambda >0$, is also a martingale.
Define the conjugate probability measure $\mathbf{P}_\lambda $ on $(\Omega ,\mathcal{F})$:
\begin{equation}
d\mathbf{P}_\lambda =M_{T\wedge n}(\lambda )d\mathbf{P}.  \label{chmeasure3}
\end{equation}
Denote $\mathbf{E}_{\lambda}$ the expectation with
respect to $\mathbf{P}_{\lambda}$.
Using the change of probability measure (\ref{chmeasure3}), we have  for all $x>0$,
\begin{eqnarray}
  \mathbf{P}\bigg( \max_{1\leq k \leq n} \frac{S_k}{V_n(\beta)}  \geq x \bigg)
 &=& \mathbf{E}_{\lambda} \left[ M_{T\wedge n}(\lambda)^{-1}\textbf{1}_{\big\{  \max_{1\leq k \leq n} S_k/V_n(\beta) \, \geq x    \big\}} \right] \nonumber \\
 &=& \sum_{k=1}^{n}\mathbf{E}_{\lambda} \bigg[\exp\bigg\{-\lambda \sum_{i=1}^k\eta_i + \Psi_{k}(\lambda) \bigg\} \textbf{1}_{\{T(x) =k\}} \bigg], \label{ghna2}
\end{eqnarray}
where
\[
\Psi_{k}(\lambda)= \sum_{i=1}^k \log \mathbf{E} \left[\exp\left\{\lambda \eta_i    \right\} \Big|\mathcal{F}_{i-1} \right].
\]
Since $\big( \eta_i,\mathcal{F}_i\big)_{i=1,...,n}$ is  conditionally symmetric, one has
$$\mathbf{E} \left[\exp\left\{\lambda \eta_i    \right\} \big|\mathcal{F}_{i-1} \right] = \mathbf{E} \left[\exp\left\{-\lambda \eta_i    \right\} \big|\mathcal{F}_{i-1} \right],$$
and thus it holds
\begin{eqnarray}\label{fsf3}
\mathbf{E} \left[\exp\left\{\lambda \eta_i    \right\} \big|\mathcal{F}_{i-1} \right] = \mathbf{E} \left[\cosh(\lambda \eta_i) \big|\mathcal{F}_{i-1} \right].
\end{eqnarray}
Since $$\cosh(x)=\sum_{k=0}^{\infty} \frac{1}{(2k)!} x^{2k}$$
is an even function, then $\cosh(\lambda \eta_i)$  is $\mathcal{F}_{i-1}-$measurable.
Thus (\ref{fsf3}) implies that $$\mathbf{E} \big[\exp\left\{\lambda \eta_i    \right\} \big|\mathcal{F}_{i-1} \big] = \cosh(\lambda \eta_i) .$$
Notice that the function $g(x)= \log \big(\cosh( x  ) \big)$ is even and convex in $x \in \mathbf{R}$ and increasing in $x \in [0, \infty)$.
Since $|\sum_{i=1}^n \eta_i| \leq n^{1- 1/\beta  } (\sum_{i=1}^n |\eta_i|^{\beta})^{1/\beta}= n^{1- 1/\beta  },$   it holds
\[
\Psi_{k}(\lambda) \leq \Psi_{n}(\lambda) = \sum_{i=1}^n g \Big(  \lambda \eta_i  \Big)  \leq n g \Big( \frac1n \sum_{i=1}^n\lambda \eta_i  \Big) \leq n g \Big( \frac{\lambda}{n^{ 1/\beta} }   \Big) .
\]
By   the fact $\sum_{i=1}^k\eta_i \geq x$   on the set $\{T(x) =k\}$, inequality (\ref{ghna2}) implies that
for all $x>0,$
\begin{eqnarray}
   \mathbf{P}\bigg( \max_{1\leq k \leq n} \frac{S_k}{V_n(\beta)}  \geq x \bigg) &\leq& \sum_{k=1}^{n}\mathbf{E}_{\lambda} \left[\exp \left\{-\lambda x +n g \Big(\frac{\lambda}{n^{ 1/\beta} }  \Big)  \right\} \textbf{1}_{\{T =k\}} \right]\nonumber\\
 &\leq&  \exp \left\{-\lambda x   +  n g \Big( \frac{\lambda}{n^{ 1/\beta} }   \Big) \right\}  \label{lastie2}.
\end{eqnarray}
The last inequality attains its minimum at
\[
\lambda=  \lambda(x)= \frac{ n^{1/\beta} }{2} \log  \bigg( \frac{n^{\frac{\beta-1}{\beta}  }+ x }{n^{\frac{\beta-1}{\beta}  }- x}\bigg) , \ \ \   \ \ \ x \in (0, n^{\frac{\beta-1}{\beta}  }).
\]
Substituting $\lambda = \lambda(x)$ in (\ref{lastie2}), we obtain the desired inequality (\ref{ineqa1}).

Notice that the function $h(x)= g(\sqrt{x})$ is  convex  and increasing in $x \in [0, \infty)$. Therefor $g(\sqrt{x})/x$ is increasing in $x,$ and
$g(\sqrt{\lambda^2/n})/(\lambda^2/n)$ is decreasing in $n.$  Thus
$$B_n(2, x)=\inf_{\lambda\geq  0}\exp \left\{-\lambda x   +  n g \Big( \frac{\lambda}{n^{ 1/2} }   \Big) \right\} $$
is increasing in $n.$
 Since $ n g ( \frac{\lambda}{n^{ 1/2} }) \rightarrow \lambda^2/2 , n\rightarrow \infty, $ we obtain
$$\lim_{n\rightarrow \infty}  B_n(2, x) = \sup_{n} B_n(2, x)    = \inf_{\lambda\geq  0}\exp \left\{-\lambda x   +  \frac{\lambda^2}{2 }   \right\}=\exp \left\{-  \frac{x^2}{2 \, }   \right\} . $$
This completes the proof of theorem.
 \hfill\qed \\

\noindent\emph{Proof of Corollary \ref{co1}.}
Since $\cosh(x) \leq \exp\{ x^2/2 \},$ we have
\begin{eqnarray}
n g \Big( \frac{\lambda}{n^{ 1/\beta} }   \Big) \ \leq \ \frac{\lambda^2}{2}  n^{1 - \frac{2}{\beta}  }
\end{eqnarray}
for all $\lambda>0.$ Thus,  from (\ref{lastie2}), for all $x > 0,$
\begin{eqnarray}
\mathbf{P}\bigg( \max_{1\leq k \leq n} \frac{S_k}{V_n(\beta)}  \geq x \bigg) &\leq& \inf_{\lambda > 0} \exp \left\{-\lambda x   +  n g \Big( \frac{\lambda}{n^{ 1/\beta} }   \Big) \right\}  \nonumber \\ &\leq& \inf_{\lambda > 0} \exp \left\{-\lambda x   + \frac{\lambda^2}{2}  n^{1 - \frac{2}{\beta}  } \right\} \nonumber \\
&=& \exp \left\{-  \frac{\ x^2}{2 }\,  n^{\frac{2}{\beta}- 1   }   \right\},
\end{eqnarray}
which gives the desired inequality (\ref{ineqa2}).\hfill\qed

\subsection*{Acknowledgements}
 The work has been partially supported by the
 National Natural Science Foundation of China (Grants no.\ 11601375).

\end{document}